\documentclass{article}
\usepackage{amssymb}

\usepackage{graphicx}
\usepackage{amsmath}

\input{tcilatex}

\begin{document}

\begin{center}
{\Huge On the number of extreme measures with fixed marginals}

\bigskip {\large M.G.Nadkarni and K.Gowri Navada}

\bigskip
\end{center}

{\large {\ }} {\large {\ \textbf{Introduction:} }}

{\large \bigskip }

{\large In his paper [3], K. R. Parthasarathy gives a bound for the number
of extreme points of the convex set of all $G-$\ invariant probability
measures on $X\times Y$\ with given marginals of full support. The purpose
of this paper is to improve this bound. }

{\large \bigskip }

{\large \textbf{Section 1:} }

{\large \bigskip }

{\large Let $X$\ and $Y$\ be finite sets with $|X|=m$\ and $|Y|=n.$\ Let $G$%
\ be a group acting on $X$\ and $Y.$\ Let $G$\ act on $X\times Y$\ \ by $%
g(x,y)=(g(x),g(y))$\ for all $g\in G$\ and $(x,y)\in X\times Y.$\ Let $X/G$\
be the set of $G$\ orbits of $X$. Write $|X/G|=m_{1},$\ $|Y/G|= n_{1}$\ and $%
|(X\times Y)/G|= m_{12}$\ . Let $\pi _{1}$\ and $\pi _{2}$\ denote the
projection maps from $X\times Y$\ to $X$\ and $Y$\ respectively. The sets $%
G(x),G(y)$\ and $G(x,y)$\ respectively denote the $G-$orbits of $x\in X,y\in
Y$\ and $(x,y)\in X\times Y.$ }

{\large Let $\mu _{1}$and $\mu _{2}$\ be $G-$\ invariant probability
measures with full support on $X$\ and $Y$\ respectively. Then $K$\ $(\mu
_{1},\mu _{2})$\ denotes the convex set of all $G$\ -invariant probability
measures $\mu $\ on $X\times Y$\ with marginals $\mu _{1}$and $\mu _{2}.$\
Note that for any measure $\mu $\ $\in $\ $K$\ $(\mu _{1},\mu _{2}),$\ the
support $S(\mu )$\ of $\mu $\ is $G-$\ invariant. Let $E(\mu _{1},\mu _{2})$%
\ denote the set of extreme points of $K$\ $(\mu _{1},\mu _{2})$. In [3],
K.R.Parthasarathy gives an estimate for the number of points in $E(\mu
_{1},\mu _{2}):$\
\begin{equation*}
|\mathbf{E}(\mu _{1},\mu _{2})|\leq \sum_{\max (m_{1},n_{1})\leq r\leq
m_{1}+n_{1}}\binom{m_{12}}{r}
\end{equation*}
$.$\ \ In this note we prove that }

{\large
\begin{equation}
|\mathbf{E}(\mu _{1},\mu _{2})|\leq \binom{m_{12}}{m_{1}+n_{1}-1}  \tag{1}
\end{equation}
}

{\large \ which considerably improves the \ above bound. Indeed $\binom{%
m_{12}}{m_{1}+n_{1}-1}$\ is one of the terms in the above sum. Moreover, if $%
G$ acts trivially or if number of $G$ orbits in $G(x)\times G(y)$ is
independent of $x$ and $\ y$, then }

{\large
\begin{equation*}
|E(\mu _{1},\mu _{2})|/\binom{m_{12}}{m_{1}+n_{1}-1}\longrightarrow 0\text{
\ \ as }m_{1},n_{1}\longrightarrow \infty .
\end{equation*}
}

{\large \bigskip }

{\large In [3], K.R. Parthasarathy has proved the following theorem: }

{\large \bigskip }

{\large \textbf{Theorem}\ ([3], Theorem 3.5): \emph{A probability measure }$%
w\in K\emph{\ }(\mu _{1},\mu _{2})$\emph{\ is extreme if and only if there
is no nonzero real valued function }$\zeta $\emph{\ on }$S(w)$\emph{\ such
that} }

{\large \emph{(i) }$\zeta (g(x),g(y))=$\emph{\ }$\zeta (x,y)$\emph{\ \ for
all }$(x,y)\in S(w),g\in G;$ }

{\large \emph{(ii) }$\sum_{y}$\emph{\ }$\zeta (x,y)$\emph{\ }$w(x,y)=0$\emph{%
\ for all }$x;$ }

{\large \emph{(iii) }$\sum_{x}$\emph{\ }$\zeta (x,y)$\emph{\ }$w(x,y)=0$%
\emph{\ for all }$y.$ }

{\large \emph{\bigskip } }

{\large \textbf{Definition :}\ A $G-$\ invariant subset $S\subset X\times Y$%
\ is said to be $G-$\emph{\ good}\ if any $G-$\ invariant real (or complex)
valued function $f$\ defined on $S$\ can be written as $f(x,y)=u(x)+v(y)$\
for all $(x,y)\in S$\ \ for some $G-$invariant functions $u$\ and $v$\ on $X$%
\ and $Y$\ respectively. }

{\large \bigskip }

{\large \textbf{Proposition 1:}\ \emph{The support }$S=\ S(\mu )$\emph{\ of
a measure }$\mu $\emph{\ }$\in K(\mu _{1},\mu _{2})$\emph{\ is }$G-$\emph{%
good if and only if }$\mu \in $\emph{\ }$E$\emph{\ }$(\mu _{1},\mu _{2}).$ }

{\large \medskip }

{\large \textbf{Proof:}\ Let $\mu \in $\ $E$\ $(\mu _{1},\mu _{2})$\ and
assume that $S$\ is not $G-$\ good. Then there exists a $G-$\ invariant
function $f$\ on $S$\ which cannot be written as $f=u+v$\ where $u$\ and $v$%
\ are $G-$\ invariant. Let $L_{G}^{2}(S,\mu )$\ denote the Hilbert space of
all $G-$invariant functions defined on $S.$\ Let $\Lambda \subset
L_{G}^{2}(S,\mu )$\ denote the set of all $G-$\ invariant functions $f$\
which have \ representation $f=u+v$\ with $u,v$\ $G-$invariant functions on $%
X$\ and $Y$\ respectively. Then $\Lambda $\ \ is a proper subspace of $%
L_{G}^{2}(S,\mu )$. Hence there exists a nonzero $\zeta \in L_{G}^{2}(S,\mu
) $\ which is orthogonal to $\Lambda .$\ Then
\begin{equation*}
\sum \zeta (x,y)u(x)\mu (x,y)=0\ and\ \sum \zeta (x,y)v(y)\mu (x,y)=0
\end{equation*}
\ for all $G-$invariant $u$\ and $v$\ defined\ on $X$\ and $Y$\
respectively. In particular
\begin{equation*}
\sum_{X}u(x)\sum_{Y}\zeta (x,y)\mu (x,y)=0\
\end{equation*}
for all $G-$invariant $u$\ on $X.$\ For any $x_{0}\in X,$\ taking $u(x)=1$\
on $G(x_{0})$\ and $u(x)=0$\ on all other orbits in $X,$\ we get
\begin{equation*}
\sum_{Y}\zeta (x_{0},y)\mu (x_{0},y)=0.\
\end{equation*}
Hence
\begin{equation*}
\ \sum_{Y}\zeta (x,y)\mu (x,y)=0
\end{equation*}
\ for all $x$\ and similarly, \
\begin{equation*}
\sum_{X}\zeta (x,y)\mu (x,y)=0\
\end{equation*}
for all $\ y$\ which contradicts the theorem above. Conversely, suppose $S$\
is $G-$good. Then $\Lambda =L_{G}^{2}(S,\mu ).$\ Let $\zeta \in
L_{G}^{2}(S,\mu )$\ satisfy the conditions of the above theorem with respect
to the measure $\mu .$\ Let $f\in L_{G}^{2}(S,\mu )$ be any function. Then $%
f $ can be written as \ $f=u+v$, for some $G-$invariant functions $u$\ and $%
v $\ on $X$\ and $Y$\ respectively.\ By condition (ii) of the above theorem,
\begin{equation*}
\sum_{X\times Y}u(x)\zeta (x,y)\mu (x,y)=\sum_{x}u(x)\sum_{y}\zeta (x,y)\mu
(x,y)=0.
\end{equation*}
}

{\large Similarly by (iii),
\begin{equation*}
\sum \zeta (x,y)v(y)\mu (x,y)\ =0.
\end{equation*}
\ Both these equations together imply
\begin{equation*}
\sum \zeta (x,y)f(x,y)\mu (x,y)=0.\
\end{equation*}
Since $f$\ is arbitrary in $L_{G}^{2}(S,\mu ),$\ $\zeta =0.$\ By the above
theorem $\mu \in E(\mu _{1},\mu _{2}),$\ which proves the proposition. }

{\large \bigskip }

{\large \textbf{Remark 1:}\ For any $x\in X$\ and $y\in Y$, the $G-$%
invariant set }$G(x)\times G(y)$ {\large can be written \ as union of }$G-$%
{\large \ orbits on }$X\times Y${\large \ whose first projection is }$G(x)$%
{\large \ and second projection is $G(y).$ } {\large \
\begin{equation*}
G(x)\times G(y)=\cup \{G(z,w)|\ \pi _{1}(G(z,w)) =G(x)\ \text{and}\ \pi
_{2}(G(z,w))\ =G(y)\}.
\end{equation*}
\ This is because the orbit $G(z,w)$\ of $(z,w)\in X\times Y$\ has $\pi
_{1}(G(z,w))=G(x)$\ and $\pi _{2}(G(z,w))$\ $=G(y)$\ if and only if $%
G(z,w)\subset G(x)\times G(y).$ }

{\large \bigskip }

{\large \textbf{Remark 2:}\ If $S$\ is a $G-$\ good set \ then $(G(x)\times
G(y))\cap S$\ contains atmost one $G-$\ orbit. This is because $S$\ cannot
contain two distinct orbits with the same projections: for, if $G(z,w)$\ and
$G(a,b)$\ are two such orbits with $\pi _{1}(G(z,w))$\ $=\pi
_{1}(G(a,b))=G(x)$\ and $\pi _{2}(G(z,w))$\ $=\pi _{2}(G(a,b))=G(y)$\ then
for any $G-$\ invariant $f=u+v$\ defined on $S$\ with $f(z,w)\neq f(a,b)$\
we will have $f(z,w)=u(z)+v(w)=u(x)+v(y)$\ and similarly, $f(a,b)=u(x)+v(y),$%
\ a contradiction. \bigskip }

{\large \textbf{Section 2: }\ }

{\large \bigskip }

{\large Let $X_{1}$ and $Y_{1}$ be two finite sets with $|X_{1}|=m_{1}$\ and
$|Y_{1}|=n_{1}.$\ \ A subset $S\subset X_{1}\times Y_{1}$\ is called \emph{%
good} \ (ref.\ [1]) if every real (or complex) valued function $f$\ on $S$\
can be expressed in the form }

{\large
\begin{equation*}
f(x)=u(x)+v(x)\text{\ for all }(x,y)\in S
\end{equation*}
}

{\large Let $X_{1}=X/G$\ and $Y_{1}=Y/G.$\ Then $X_{1}\times Y_{1}$\ can be
identified with a set whose points are $G(x)\times G(y),$\ $x\in X,y\in Y.$\
The $G-$invariant functions on $X$\ and on $Y$\ are in one-to-one
correspondence with the functions on $X_{1}$\ and $Y_{1}$\ respectively. }

{\large Let $\widetilde{S}\subset (X\times Y)/G$\ denote the set of all $G-$%
orbits in a $G-$invariant subset $S$\ of $X\times Y.$\ Define $\phi :%
\widetilde{S}\longrightarrow X_{1}\times Y_{1}$\ by $\phi
(G(x,y))=G(x)\times G(y)$. }

{\large One can show that subsets of good sets are good and every good set $%
S\subset X_{1}\times Y_{1}$\ is contained in a maximal good subset of $%
X_{1}\times Y_{1}$. Further every maximal good set of $X_{1}\times Y_{1}$\
contains $m_{1}+n_{1}-1$\ elements. (ref [1])}

{\large \bigskip }

{\large \textbf{Proposition 2: }\emph{A }$G-$\emph{invariant subset }$%
S\subset X\times Y$\emph{\ is }$G-$\emph{good if and only if }$\phi $\emph{\
is one-to-one on }$\widetilde{S}$\emph{\ and }$\phi (\widetilde{S}$\emph{\ }$%
)$\emph{\ is good in }$X_{1}\times Y_{1}$\emph{.\ Further, }$S$\emph{\ is
maximal }$G-$\emph{good set if and only if }$\phi $\emph{\ is one-to-one on }%
$S$\emph{\ and }$\phi (\widetilde{S}$\emph{\ }$)$\emph{\ is maximal good set
in }$X_{1}\times Y_{1}$\emph{.} }

{\large \bigskip}

{\large \textbf{Proof:}\ Assume\ $S$\ is $G-$good. By remark 2, if $S$\ is $%
G-$\ good, then\ $\phi $\ is one-to-one on $\widetilde{S}.$\ Let $f$\ be any
real (or complex) valued function defined on $\phi (\widetilde{S}$\ $).$\
Define $g$\ on $\widetilde{S}$\ by $g=f\circ \phi .$\ This map $g$\ gives
rise \ to a $G-$invariant map on $S,$\ again denoted by $g.$\ Writing $g=u+v$%
, where $u$\ and $v$\ are\ $G-$invariant functions on $X$\ and $Y$\ \
respectively,\ and noting that $u$\ and $v$\ are constant on each orbit, we
can define $\widetilde{u}$\ and \ $\widetilde{v}$\ on $X_{1}$\ and $Y_{1}$\
\ by \ $\widetilde{u}(G(x))\ =u(x)$\ and\ \ $\widetilde{v}(G(y))\ =v(y).$\
It is easy to see that $f=\widetilde{u}+$\ $\widetilde{v}.$\ So \ $\phi (%
\widetilde{S})$\ is good. Conversely, let $S\subset X\times Y$ \ be such
that $\phi $ is one-to-one on $\widetilde{S}$ and $\phi (\widetilde{S}$\ $)$%
\ is good. Since $\phi $\ is one-to-one, any $G(x)\times G(y)$\ intersects$\
S$\ in atmost one orbit. Given a function $g$\ on $S$\ we can define $f$\ on
$\phi (\widetilde{S})$\ as $f=g\circ \phi ^{-1}.\ $Since $f$ is defined on
the good set $\phi (\widetilde{S})$ we can write $f$ as $f=\widetilde{u}+$\ $%
\widetilde{v}$ where $\widetilde{u},$\ $\widetilde{v}$ are defined on $\pi
_{1}(\phi (\widetilde{S}))$ and $\pi _{2}(\phi (\widetilde{S}))$
respectively.\ Defining $u(x)=\widetilde{u}(G(x))\ $and $v(y)=\widetilde{v}%
(G(y))\ $we get $G-$ invariant functions $u$ and $v$\ with $g=u+v.$\ \ Now
suppose $S$\ is a maximal $G-$good set. We know from the first part of the
theorem that $\phi $\ is one-to-one on $\widetilde{S}.$ \ If $\phi (%
\widetilde{S}$\ $)$\ is not a maximal good set, there exists a point, say $%
G(a)\times G(b)\notin \phi (\widetilde{S}$\ $),$\ such that $\phi (%
\widetilde{S})\cup \{G(a)\times G(b)\}$\ is good. Then, since $\{G(a)\times
G(b)\}\cap S=\emptyset ,$ the map $\phi $\ is one-to-one on $\widetilde{T}$
\ where \ $T=G(a,b)\cup S.$ \ Using the first part of the theorem, $T$\ is $%
G-$good contradicting the maximality of $S.$\ The converse can be proved in
a similar manner. This completes the proof of the proposition. }

{\large \bigskip }

{\large By corollary 3.6 of [3], different extreme points of $K(\mu _{1},\mu
_{2})$\ have distinct supports.\ As pointed out by the referee, this fact is
also a consequence of proposition 1: Assume that }$\mu ,\nu \in ${\large \ $%
E $\ $(\mu _{1},\mu _{2}),$ with }$\mu \neq \nu ,${\large \ having the same
support }$\ S.${\large \ By proposition 1 }$S${\large \ is }$G-${\large %
good. But this is a contradiction since }$S${\large \ is also the support of
}$(\mu +\nu )/2,${\large \ which is not extreme. \ Further, for $\mu $\ and $%
\nu \in $\ $E$\ $(\mu _{1},\mu _{2})$\ the measure $(\mu +\nu )/2\in K$\ $%
(\mu _{1},\mu _{2})$\ is not extreme, and so by Proposition 1 its support $%
S(\mu )\cup S(\nu )$\ is not a $G-$\ good set. Further, for $\mu $\ and $\nu
\in $\ $E$\ $(\mu _{1},\mu _{2})$\ the measure $(\mu +\nu )/2\in K$\ $(\mu
_{1},\mu _{2})$\ is not extreme, and so by Proposition 1 its support $S(\mu
)\cup S(\nu )$\ is not a $G-$\ good set. This shows that supports\ of
different measures in $E(\mu _{1},\mu _{2})$\ are contained in different
maximal $G-$good sets of $X\times Y:$ \ Because, if $\mu \neq \nu \in E\
(\mu _{1},\mu _{2})$ such that $S(\mu )$}$\subset ${\large \ }$S${\large \
and $S(\nu )$}$\subset ${\large \ }$S${\large \ for some maximal }$G-$%
{\large good set }$S$ {\large then the measure} {\large $(\mu +\nu )/2\in K$%
\ $(\mu _{1},\mu _{2})$ has its support $S(\mu )\cup S(\nu )$ contained in $%
S.$ Since $S$ is $G-$good,} {\large $S(\mu )\cup S(\nu )$ is also }$G-$%
{\large good a contradiction to proposition 1 as $(\mu +\nu )/2$ is not
extreme. Therefore, $|E(\mu _{1},\mu _{2})|$\ is bounded by the number of
maximal $G-$good sets of $X\times Y.$ }

{\large \bigskip}

{\large Let $S$ be a maximal $G-$good set in $X\times Y$. By Proposition 2, $%
\phi (\widetilde{S})$ is a maximal good set in $X_{1}\times Y_{1}$. Since $%
\phi $\ is one-to-one on $\widetilde{S}$, \ $\widetilde{S}$ contains $%
m_{1}+n_{1}-1$ orbits of $G$. Since the number of orbits in $X\times Y$\ is $%
m_{12},$ and any maximal $G-$good set in $X\times Y$ is of the form $\phi (%
\widetilde{S}$\ $),$ \ the total number of maximal $G-$good sets in $X\times
Y$\ is less than or equal to $\binom{m_{12}}{m_{1}+n_{1}-1}.$ This proves $%
(1).$ }

{\large \bigskip }

{\large We give an example to show that the above bound is sharp. Let }$G$%
{\large \ be the group }$S_{n},${\large \ the permutation group on }$n$%
{\large \ elements. Let }$\ X=\{1,2,...,n\}${\large \ and }$Y${\large \ be
the set }{\large $S_{n}$ . Here $|X|=n$}{\large \ and $|Y|=n!.$}{\large \
Then }$G${\large \ acts on }$X${\large \ in the obvious manner and on }$Y$%
{\large \ by }$g(h)=g\circ h.${\large \ The only }$G-${\large \ invariant
subset of }$X${\large \ is }$X${\large \ itself and the only }$G-${\large \
invariant subset of }$Y${\large \ is }$Y${\large \ itself. Then }$G${\large %
\ also acts on }$X\times Y${\large \ diagonally. That is, }$%
g(x,y)=(g(x),g(y)).${\large \ } {\large For any }$(x,y)\in X\times Y,$%
{\large \ the set }$G(x,y)=${\large \ }$\{(g(x),g(y))|g\in G\}${\large \ is
a }$G-${\large \ invariant subset of }$X\times Y${\large \ with }$n!${\large %
\ number of elements and }$\ G(x)\times \ G(y)${\large \ is the whole set }$%
X\times Y.${\large \ In this case, }$|X/G|=m_{1}=1${\large \ and \ }$%
|Y/G|=n_{1}=1${\large \ and }$|(X\times Y)/G|=m_{12}${\large \ }$=n.${\large %
\ Therefore, }$\binom{m_{12}}{m_{1}+n_{1}-1}=n.$

{\large The only }$G-${\large \ invariant probability measures on }$X$%
{\large \ and }$Y${\large \ are uniform measures. That is, }$\mu _{1}(x)=%
\frac{1}{n}${\large \ for all }$x\in X${\large \ and }$\mu _{2}(y)=\frac{1}{%
n!}${\large \ for all }$y\in Y.${\large \ So the only }$G-${\large \
invariant functions on }$X${\large \ and }$Y${\large \ are constant
functions.\ If }$\mu \in ${\large \ }$E(\mu _{1},\mu _{2}),${\large \ then
the support }$S${\large \ of }$\mu ${\large \ should be }$G-${\large \ good.
Any }$G-${\large \ invariant function }$f${\large \ defined on }$S,${\large %
\ can be written as }$f=u+v${\large \ where }$u,v${\large \ are }$G-${\large %
\ invariant functions on }$X${\large \ and }$Y${\large \ respectively. This
shows that }$f${\large \ must be constant, which means }$S${\large \
consists of a single orbit, say }$S=G(x,y).${\large \ Then }$\mu
((g(x),g(y))=\frac{1}{n!}${\large \ for all }$g\in G.${\large \ Observe that
the collection }$\{g(y)|g\in G\}${\large \ has all }$n!${\large \ different
elements\ whereas in the collection }$\{g(x)|g\in G\}${\large \ \ every
value of }$g(x)${\large \ is repeated }$(n-1)!${\large \ times. This shows
that every such uniform measure }$\mu ${\large \ supported on any single
orbit }$G(x,y)${\large \ has marginals }$\mu _{1}${\large \ and }$\mu _{2}.$%
{\large \ Since there are }$n${\large \ orbits in }$X\times Y,${\large \ we
get }$|E(\mu _{1},\mu _{2})|=n.$

\bigskip

{\large Now we state some results about good subsets of \ $X_{1}\times Y_{1}
$\ not necessarily $G$-good sets ( ref. [1], [2] ). }

{\large Consider any two points $(x,y),(z,w)\in S\subset X_{1}\times Y_{1}$\
where $S$\ is any (not necessarily good) subset of $X_{1}\times Y_{1}.$\ We
say that $(x,y),(z,w)$\ are\emph{\ linked} if \ there exists a sequence of
points $(x_{1},y_{1})=(x,y),(x_{2},y_{2})...(x_{n},y_{n})=(z,w)$\ of points
of $S$\ such that }

{\large $(i)$\ for any $1\leq i\leq n-1$\ exactly one of the following
equalities hold: }

{\large $x_{i}=x_{i+1}$\ or $y_{i}=y_{i+1};$ }

{\large $(ii)$\ if $x_{i}=x_{i+1}$\ then $y_{i+1}=y_{i+2},$\ and if $%
y_{i}=y_{i+1}$\ then $x_{i+1}=x_{i+2},1\leq i\leq n-2.$ }

{\large We also call this a \emph{link} joining $(x,y)$\ to $(z,w). $\ A
nontrivial link joining $(x,y)$\ to itself is called a\ \emph{loop}. }

{\large \bigskip }

{\large \textbf{Theorem}\ (ref. [1], cor. 4.11): A subset $S\subset
X_{1}\times Y_{1}$\ is good if and only if $S$\ contains no loops. }

{\large \bigskip }

{\large \bigskip }

{\large \textbf{Remark 3:} Let the orbits in $\widetilde{S}$\ be
\begin{equation*}
G(x_{1},y_{1}),G(x_{2},y_{2}),...,G(x_{m_{1}+n_{1}-1},y_{m_{1}+n_{1}-1}).
\end{equation*}
\ Then $S\cap (G(x_{i})\times G(y_{i}))$ $=G(x_{i},y_{i})$\ for $1\leq i\leq
m_{1}+n_{1}-1.$ Let $G(z,w)$\ be any other orbit in $G(x_{i})\times G(y_{i})$
and let $S^{\prime }=(S\setminus G(x_{i},y_{i}))\cup $ $G(z,w).$ It is clear
that $S^{\prime }$\ is maximal $G-$good set with $\phi (\widetilde{S})=\phi (%
\widetilde{S^{\prime }})$. If $\alpha _{i}$\ denote the number of orbits in $%
G(x_{i})\times G(y_{i}),$ then\ there are $\alpha _{1}\alpha _{2}...$\ $%
\alpha _{m_{1}+n_{1}-1}$\ many maximal $G-$good sets in $X\times Y$ with
image $\phi (\widetilde{S})$\ under $\phi .$ }

{\large It seems likely that $|E(\mu _{1},\mu _{2})|/\binom{m_{12}}{%
m_{1}+n_{1}-1}\longrightarrow 0$ \ \ as $m_{1},n_{1}\longrightarrow \infty .$
We show this in the case $G=(e)$ and more generally when number of $G$
orbits in $G(x)\times G(y)$ is independent of $x$ and $\ y.$ For that we
first prove the following theorem. }

{\large \bigskip }

{\large \textbf{Theorem:}\ \emph{Let }$X$\emph{\ }$=\{x_{1},x_{2},...x_{m}\}
$\emph{\ and }$Y=\{y_{1},y_{2},...y_{n}\}$\emph{\ be two finite sets. Then} }

{\large $(i)$\emph{\ the number of maximal good sets contained in }$X\times
Y $\emph{\ is }$m^{n-1}n^{m-1}.$ }

{\large $(ii)$\emph{\ the number of maximal good sets among them with
exactly }$k$\emph{\ fixed points, (say }$%
(x_{i},y_{j_{1}}),...(x_{i},y_{j_{k}})$\emph{\ ) having a fixed first
coordinate say }$x_{i}$\emph{\ is : }$kn^{m-2}(m-1)^{n-k},$\emph{\ }$1\leq
k\leq n.$ }

{\large $(iii)$\emph{\ the number of maximal good sets with exactly }$k$%
\emph{\ fixed points having a fixed second coordinate say }$y_{j}$\emph{\
is: }$km^{n-2}(n-1)^{m-k},$\emph{\ }$1\leq k\leq m.$ }

{\large \bigskip }

{\large \textbf{Proof:}\ \ We use induction on $m+n.$\ The result is true
for $m=1$\ and $n=1$. Assume the result for all values of $|X|\leq m$\ and $%
|Y|\leq n.$\ We prove the result for $|X|=m$\ and $|Y|=n+1.$\ Let $X$\ $%
=\{x_{1},x_{2},...x_{m}\}$\ and $Y=\{y_{1},y_{2},...y_{n},y_{n+1}\}. $\
Consider a $m\times (n+1)$\ grid of $\ m(n+1)$\ cells with $m$\ rows
corresponding to $\{x_{1},x_{2},...,x_{m}\}$\ and $n+1$\ columns
correspondiong to $\{y_{1},y_{2},...,y_{n+1}\}$\ Associate $(i,j)$th cell
with the point $(x_{i},y_{j})\in X\times Y.$\ We say that $(x_{i},y_{j})\in
(i,j)$th cell. }

{\large To prove $(iii)$\ let $S$\ be a maximal good set in $X\times Y.$\
Then $|S|=m+n.$\ Suppose $S$\ contains exactly $k$\ points with fixed second
coordinate, say $y_{n+1}$. Without loss of generality we assume them to be $%
(x_{1},y_{n+1}),$\ $(x_{2},y_{n+1}),$\ $...$\ $(x_{k},y_{n+1}).$\ Denote
\begin{equation*}
K=\{(x_{1},y_{n+1}),(x_{2},y_{n+1}),...(x_{k},y_{n+1})\}.\
\end{equation*}
}

{\large (i) Atleast one of these first $k$\ rows contain atleast two points
of $S,$\ i.e., there exist a point $(x_{i},y_{j})$\ of $S$\ with $1\leq
i\leq k$\ and $1\leq j\leq n.$ }

{\large \emph{Proof: }Otherwise leaving these $k$\ rows and the last column,
the remaining points of $S$\ will be a good set with $m+n-k$\ points using $%
m+n-k$\ coordinates which is not possible. }

{\large (ii) \ If $(x_{i},y_{j})\in S$\ with $1\leq i\leq k$\ and $1\leq
j\leq n.$\ Then the $j$the column (which contains the point $(x_{i},y_{j})$)
has no other point $(x_{l},y_{j})$\ of $S$\ with $1\leq l\neq i\leq k$\
because the four points $%
\{(x_{i},y_{j}),(x_{i},y_{n+1}),(x_{l},y_{n+1}),(x_{l},y_{j})\}$\ form a
loop. }

{\large (iii) Suppose $(x_{i},y_{j})\in S$\ \ for some $1\leq i\leq k$\ and $%
1\leq j\leq n.$\ Then the set got by dropping the point $(x_{i},y_{j})$\ and
adding $(x_{l},y_{j}),1\leq l\neq i\leq k$\ to $S$\ clearly contain no loop
and so is maximal good. }

{\large Let $S^{\prime }$\ be the maximal good set obtained in this way by
replacing all the points $(x_{i},y_{j}),$\ $1\leq i\leq k$\ and $1\leq j\leq
n$\ of $S$\ \ by $(x_{1},y_{j}),1\leq j\leq n$\ . }

{\large Then each of the rows corresponding to $x_{2},...,x_{k}$\ \ contains
exactly one point of $S^{\prime }.$\ The set $S^{\prime \prime }$\ got from $%
S^{\prime }$\ by dropping these rows and the last column will be a maximal
good set in $\{x_{1},x_{k+1},...,x_{m}\}\times \{y_{1},y_{2},...y_{n}\}$\
and contains $m+n-k$\ elements. }

{\large $\ $\ By induction hypothesis, the number of maximal good sets in $%
\{x_{1},x_{k+1},...,x_{m}\}\times \{y_{1},y_{2},...y_{n}\}$\ having exactly $%
r$\ points in $r$\ fixed positions in the first row, is: $%
rn^{m-k-1}(m-k)^{n-r},$\ for $1\leq r\leq n.$\ Consider any such maximal
good set, say $A.$\ Further add the dropped rows and the last column.
Enlarge $A$\ by adding the first $k$\ points of the $(n+1)$th column, call
this set $B$. It is a maximal good set in $\{x_{1},x_{2},...,x_{m}\}\times
\{y_{1},y_{2},...y_{n+1}\}$. Any point $(x_{1},y_{j}),1\leq j\leq n$\ in $B$%
\ can be replaced by $(x_{l},y_{j}),$\ for any $1\leq l\leq k$\ and the
resulting set will continue to remain maximal good in $%
\{x_{1},x_{2},...,x_{m}\}\times \{y_{1},y_{2},...y_{n+1}\}$. In this way
each one of $rn^{m-k-1}(m-k)^{n-r}$\ maximal good set $A$\ gives rise to $%
k^{r}$\ maximal good sets in the original $m\times (n+1)$\ matrix. .
Further, we can choose the $r$\ points in the first row in $\binom{n}{r}$\
ways . Adding over $r$, the total number of maximal good sets with exactly $%
k $\ cells in $k$\ fixed positions of the last column \ is: $\ $%
\begin{equation*}
\sum_{r=1}^{n}\binom{n}{r}k^{r}rn^{m-k-1}(m-k)^{n-r}=kn^{m-k-1}n%
\sum_{r=0}^{n-1}\binom{n-1}{r}k^{r-1}(m-k)^{n-r}
\end{equation*}
}

{\large
\begin{equation*}
=kn^{m-k}(m-k+k)^{n-1}=kn^{m-k}m^{n-1}
\end{equation*}
which is $(iii)$\ \ for $m\times (n+1)$\ matrix. }

{\large To prove $(i),$\ since we can choose the $k$\ points in the last
column in $\binom{m}{k}$\ ways, the total number of maximal good sets with
exactly $k$\ points from the last column is: $\binom{m}{k}$\ $%
kn^{m-k}m^{n-1}.$\ \ The total number of maximal good sets in \ $X\times Y$\
is got by adding these numbers as $k$\ varies from $1$\ to $m:$\ $\ \ \ \ $%
\begin{equation*}
\sum_{k=1}^{m}\binom{m}{k}\ kn^{m-k}m^{n-1}\ =\ \sum_{k=0}^{m-1}\binom{m-1}{k%
}\ n^{m-k-1}m^{n}\ =\ \ m^{n}(n+1)^{m-1}.\ \ \
\end{equation*}
$\ $ }

{\large $(ii)$\ can be proved in a similar way as $(iii)$. $\ $This
completes the proof of the theorem. }

\bigskip

{\large Next, we prove that as $m,n\rightarrow \infty $\ the ratio $%
m^{n-1}n^{m-1}/\binom{mn}{m+n-1}\rightarrow 0$\ as $m,n\rightarrow \infty .$
\
\begin{equation*}
\lim_{m,n\rightarrow \infty }m^{n-1}n^{m-1}/\binom{mn}{m+n-1}=
\end{equation*}
}

\begin{equation*}
\lim_{m,n\rightarrow \infty }m^{n-1}n^{m-1}(m+n-1)!((mn-m-n+1)!/(mn)!\
\end{equation*}

{\large By Sterling's formula, we know that $n!\thicksim \frac{\sqrt{2\pi }%
n^{n+\frac{1}{2}}}{e^{n}}$ for large $n.$ \bigskip\ Using this expression
one can show that }

{\large
\begin{equation}
m^{n-1}n^{m-1}/\binom{mn}{m+n-1}\leq C\frac{(1-\frac{1}{m})^{mn}(1+\frac{n}{m%
})^{m}(1-\frac{1}{n})^{mn}(1+\frac{m}{n})^{n}}{(1-\frac{1}{m})^{n}(1-\frac{1%
}{n})^{m}(m+n)^{\frac{1}{2}}}  \tag{2}
\end{equation}
}

{\large for some constant $C.$ If $\frac{m}{n}\geq 1$, since $(1-\frac{1}{m}%
)^{m}$\ increases to $e^{-1}$, the right hand side of $(2)$\ tends to $0$\
as $m,n\longrightarrow \infty .$ \ }

{\large The case where $\frac{m}{n}$$\leq 1$ is similar because the the
expression on the right hand side of $(2)$ is symmetric with respect to $m$
and $n.$ }

{\large If $G=(e),$ the maximal $G$-good sets in $X\times Y$ are just the
maximal good sets and the number of maximal good sets, by the previous
theorem\ is, $m^{n-1}n^{m-1}.$ In this case $m_{12}=m_{1}n_{1}.$ Therefore
\begin{equation*}
|E(\mu _{1},\mu _{2})|/\binom{m_{12}}{m_{1}+n_{1}-1}\leq m^{n-1}n^{m-1}/%
\binom{m_{1}n_{1}}{m_{1}+n_{1}-1}\longrightarrow 0
\end{equation*}
}

{\large as $m$ and $n\longrightarrow \infty .$ }

{\large Now suppose that the number of $G$-orbits in $G(x)\times G(y)$ is a
constant, say $a,$ for all $x$ and $y.$ Then by remark 3, the number of
maximal $G$-good sets in $X\times Y$ is $%
a^{m_{1}+n_{1}-1}m_{1}^{n_{1}-1}n_{1}^{m_{1}-1}$ and $m_{12}=am_{1}n_{1}.$
Therefore,
\begin{equation*}
|E(\mu _{1},\mu _{2})|/\binom{m_{12}}{m_{1}+n_{1}-1}\leq \left(
a^{m_{1}+n_{1}-1}m_{1}^{n_{1}-1}n_{1}^{m_{1}-1}\right) /\binom{am_{1}n_{1}}{%
m_{1}+n_{1}-1}
\end{equation*}
\begin{equation*}
\leq \left( a^{m_{1}+n_{1}-1}m_{1}^{n_{1}-1}n_{1}^{m_{1}-1}\right)
/a^{m_{1}+n_{1}-1}\binom{m_{1}n_{1}}{m_{1}+n_{1}-1}
\end{equation*}
\ }

{\large
\begin{equation*}
=\left( m_{1}^{n_{1}-1}n_{1}^{m_{1}-1}\right) /\binom{m_{1}n_{1}}{%
m_{1}+n_{1}-1}\longrightarrow 0
\end{equation*}
}

{\large as $m$ and $n\longrightarrow \infty .$ }

{\large \bigskip }

{\large \textbf{Note:}\ The maximal good sets in $X\times Y$\ can be
associated in a one-to-one manner with the spanning trees of a complete
bipartite graph. Consider the complete bipartite graph $K_{m,n}$\ where $%
|X|=m$\ and $|Y|=n.\ $\ A subset $S\subset X\times Y\ $\ is maximal good if
and only if $\ |S|=m+n-1$\ and in the grid corresponding to $X\times Y$, $S$%
\ contains no loops. Construct an $m\times n$\ matrix corresponding to any
spanning tree $T$\ in $\ \ K_{m,n}$\ as follows: Identifying the elements of
$X$\ and $Y$\ with the veritces of $K_{m,n},$\ let $V=(X,Y)$\ denote the
vertices of $K_{m,n}.$\ Whenever the edge $(x_{i},y_{j})\in T,$\ put $(i,j)$%
th entry in the matrix equal to one$;$\ otherwise $(i,j)$th entry is zero.
Since $T$\ is a spanning tree, there are exactly $m+n-1$\ nonzero entries in
the matrix. As $T$\ contains no cycles, the nonzero entries in the matrix
donot form a loop. Therefore the nonzero entries of the matrix correspond to
a maximal good set in the grid corresponding to $X\times Y.$\ This
correspondence is one-to-one. In [5], it is proved that the number of
spanning trees of $K_{m,n}$\ is $m^{n-1}n^{m-1}.$\ But the proof makes use
of the determinant of the matrix and is different from the one given here. }

\bigskip

{\large \textbf{Acknowledgement:} The second author thanks CSIR for funding
the project of which this paper forms a part. Sincere thanks to IMSc,
Chennai and ISI, Bangalore for providing short visiting appointments.}

{\large \bigskip }

{\large \textbf{References:} }

{\large \medskip }

{\large [1]\quad Cowsik R C, Klopotowski A and Nadkarni M G, When is $%
f(x,y)=u(x)+v(y)$?, Proc. Indian Acad. Sci. (Math. Sci.) 109 (1999) 57-64. }

{\large [2] \ \ Nadkarni M G, Kolmogorov's superposition theorem and sums of
algebras, The Journal of Analysis vol. 12 (2004) 21-67. }

{\large [3]\quad K. R Parthasarathy, Extreme points of the convex set of
joint probability distributions with fixed marginals, Proc. Indian Acad.
Sci. (Math. Sci.) 117 (2007) 505-516. }

{\large [4] \ \ \ D. B. West, Introduction to Graph Theory, Second Edition,
Pearson Education. }

{\large [5] \ \ \ Lovasz L., Combinatorial problems and exercises, Second
Edition, Elsevier Science Ltd. }

{\large \bigskip }

{\large \textbf{Address:} }

{\large \medskip }

{\large 1. Prof. M.G. Nadkarni }

{\large Department of Mathematics }

{\large University of Mumbai }

{\large Kalina Campus, Santacruz East }

{\large Mumbai-400098 }

{\large email: mgnadkarni@gmail.com }

{\large \bigskip }

{\large 2. K. Gowri Navada }

{\large Department of Mathematics }

{\large Periyar University }

{\large Salem-636011 }

{\large Tamil Nadu }

{\large email: gnavada@yahoo.com }

{\large \bigskip }

\end{document}